\documentclass{amsart}

\usepackage{amssymb}

\newtheorem{defi}{Definition}
\newtheorem{lemma}{Lemma}
\newtheorem{theorem}{Theorem}
\newtheorem{proposition}{Proposition}

\usepackage{amsmath}
\usepackage{amsfonts}
\usepackage{amssymb}

\DeclareMathOperator{\ord}{ord}

\newcommand{\Le}{\leqslant}
\newcommand{\Ge}{\geqslant}

\newcommand{\K}{{\mathbf k}}

\newcommand{\N}{{\mathbb N}}

\newcommand{\Q}{{\mathbb Q}}
\newcommand{\C}{{\mathcal C}}

\def\B{{\mathcal B}}

\newcommand{\ld}{\mathop{\rm ld}\nolimits}

\newcommand{\rk}{\mathop{\rm rk}\nolimits}

\makeatletter
\def\Ddots{\mathinner{\mkern1mu\raise\p@
\vbox{\kern7\p@\hbox{.}}\mkern2mu
\raise4\p@\hbox{.}\mkern2mu\raise7\p@\hbox{.}\mkern1mu}}
\makeatother

\begin{document}

\title{On the generalised Ritt problem as a computational problem}
\thanks{The work was partially supported by the Russian Foundation for Basic Research, project no. 05-01-00671.}
\author{Oleg Golubitsky}
\address{University of Western Ontario, Department of Computer Science, London, Ontario, Canada, N6A 5B7}
\email{Oleg.Golubitsky@gmail.com}
\thanks{The first author was also partially supported by NSERC Grant PDF-301108-2004.}
\author{Marina Kondratieva}
\address{Moscow State University, Department of Mechanics and Mathematics, Leninskie gory, Moscow, Russia, 119991}
\email{kondratieva@sumail.ru}
\author{Alexey Ovchinnikov}
\address{University of Illinois at Chicago, Department of Mathematics, Statistics, and Computer Science, Chicago, IL 60607-7045, USA}
\email{aiovchin@math.uic.edu}
\thanks{The third author was also partially supported by NSF Grant CCR-0096842.}

\dedicatory{To our advisor, Eugeny V. Pankratiev}


\begin{abstract} The Ritt problem asks if there is an algorithm that tells 
whether one prime differential ideal is contained in another one if both are 
given by their characteristic sets. We give several equivalent formulations 
of this problem. In particular, we show that it is equivalent to testing 
if a differential polynomial is a zero divisor modulo a radical differential ideal. 
The technique used in the proof of equivalence yields algorithms for computing a 
canonical decomposition of a radical differential ideal into prime components and
a canonical generating set of a radical differential ideal. Both proposed 
representations of a radical differential ideal are independent of the given 
set of generators and can be made independent of the ranking.
\end{abstract}

\maketitle

\section{Introduction}
The Ritt problem is an algebraic treatment of the following phenomenon
occurring with differential equations. One can show that the solution set
to the differential equation $y'^2-4y=0$ consists of a family of parabolae
and the zero special solution. The solution set to the equation $y'^2-4y^3 =0$
consists of hyperbolae and the zero solution as well. The major difference between
these two solution sets is that in the first case the zero solution is not a limit
of the family while in the second case it is. In the algebraic language, the first
equation does not generate a prime differential ideal, while the second one
does. The Ritt problem is a generalisation of this phenomenon to the case
of several non-linear algebraic PDEs.

 More precisely, consider a differential polynomial ring over a differential
  field of characteristic zero.
  J.F.~Ritt has shown that every radical differential ideal
  can be uniquely represented as a finite intersection of its 
  essential prime components, which are the minimal prime 
  differential ideals containing it.
  We call the problem of computing the characteristic sets of 
  the essential prime components, given a finite set of
  generators of a radical differential idea, the generalised
  Ritt problem. 
  As of today, the problem remains largely unsolved, with some
  important special cases having been treated by J.F.~Ritt, 
  H.~Levi, E.R.~Kolchin, and R.M.~Cohn \cite{Kol,Levi,Rit,Cohn1,Cohn2}. 

  We give a precise formulation of this problem as a 
  computational problem, emphasise its dependence on the 
  differential field, and then reformulate the problem 
  to make it field-independent. One of the versions of the Ritt problem
  that we discuss here is given in terms of zero-divisors modulo an
  ideal. In particular, when we test primality by definition for a
  pair of polynomials each not in the ideal we need to determine whether
  their product is in the ideal. This sentence has two universal
  quantifiers, one for each polynomial. We show that it
  is sufficient to quantify with respect to just one polynomial. 
  
  To illustrate these ideas
  we also exhibit an algorithm whose input is a finite set $F$ of differential polynomials
  and the output is a collection of prime differential ideals given by their characteristic
  sets such that:
  \begin{itemize}
  \item the radical differential ideal $I$ generated by $F$ is the intersection of the prime differential ideals and
  \item this set of prime differential ideals is uniquely determined just by the ideal $I$
  and does not depend on the chosen set of generators $F$.
  \end{itemize}   We hope that the
  approach we propose in this paper will bring us closer to determining decidability of the Ritt problem.
  
\section{Basic definitions}
Let $\K$ be a differential field of characteristic zero 
with derivations $\Delta=\{\delta_1,\ldots,\delta_m\}$.
Let $Y=\{y_1,\ldots,y_n\}$ be a finite set of differential
indeterminates and $$\Theta Y := \left\{\delta_1^{i_1}\cdot\ldots\cdot\delta_m^{i_m}y_j\:|\:i_k\Ge 0, 1\Le k\Le m;\: 1\Le j\Le n\right\}.$$ The ring
of differential polynomials $\K\{Y\}$ is the ring of commutative
polynomials $\K[\Theta Y]$ with the natural structure of a $\Delta$-ring.

Let $I$ be a radical differential ideal in $\K\{Y\}$, that is,
a radical ideal stable under the action of $\Delta$.
A prime differential ideal $P$ is called an 
{\it essential prime component}
of $I$ if $I\subseteq P$ and there is no prime differential
ideal $Q$ satisfying $I\subseteq Q\subsetneq P$.
According to \cite{Rit}, every radical differential ideal $I$
has finitely many essential prime components and equals
their intersection; every prime differential ideal containing $I$
must contain an essential prime component of $I$; and every
set of prime differential ideals whose intersection equals $I$
contains the set of essential prime components of $I$ as
a subset. The {\it generalised Ritt problem} is: given a 
radical differential ideal, find its essential prime components.

To formulate it as a computational problem, we need to 
fix an encoding of the input (a radical differential ideal in 
$\K\{Y\}$), and of the output (a finite set of prime differential 
ideals in $\K\{Y\}$).
A radical differential ideal is usually given by a finite
generating set of differential polynomials. To be able
to compute with them, we have to assume that $\K$ is a computable 
differential field. The following definition is directly derived 
from \cite[Definition 5]{Rabin}, see \cite{Miller08}:

\begin{defi} 
   A differential ring $R$ is called computable, if 
   there exists an embedding 
   $$i: R \to \N$$
   such that
   \begin{enumerate}
      \item The image $i(R)$ is a recursive subset of $\N$
      \item The addition, multiplication, and derivation
        functions induced on this image by the addition,
        multiplication, and derivations on $R$ are 
        computable functions.
   \end{enumerate}
\end{defi}

The embedding $i$ in the above definition is called
an {\it admissible indexing} of $R$.
If $\K$ is a computable differential field, then the ring of 
differential polynomials $\K\{Y\}$ is a computable differential 
ring: the value of an admissible indexing $j:\K\{Y\}\to \N$ can 
be effectively computed for each differential polynomial 
$f\in\K\{Y\}$, given the indices of its coefficients w.r.t. 
the indexing $i:\K\to\N$ \cite{Rabin}. 

Having fixed such an indexing $j$, we can unambiguously
say that a differential polynomial $f\in\K\{Y\}$ 
is {\it given}, meaning that actually $j(f)$ is given.
A prime differential ideal can be given either by 
a generating set (since it is radical), or by a 
characteristic set. The latter is defined as follows.
Define a ranking on
$\Theta Y$
as a total
order relation $\le$ satisfying:
$$u\le\delta u,\;\;\;u\le v\;\Rightarrow\delta u\le\delta v$$
for all $u,v\in\Theta Y$, $\delta\in\Delta$.
For a differential polynomial $f\in\K\{Y\}\setminus\K$, 
call the derivative of the highest rank effectively present 
in $f$ the leader of $f$, denoted $u=\ld_\le f$.
If the degree of $f$ w.r.t. $u$ is $d$, then $u^d$ is
called the rank of $f$, denoted $\rk_\le f$.
The ranking then extends to the set of ranks:
$$u_1^{d_1}\le u_2^{d_2}\;\iff\;u_1<u_2\;{\rm or}\;(u_1=u_2\;{\rm and}\:d_1\le d_2)$$
and to the set of all finite sets of ranks:
$R_1\le R_2$ iff either $R_1=R_2$, or else the minimal element 
of the symmetric difference 
$(R_1\setminus R_2)\cup (R_2\setminus R_1)$ belongs to $R_1$.
All three relations: $\le$ on the set of derivatives, 
$\le$ on the set of ranks, and $\le$ on the set of finite sets 
of ranks, are well-orderings.

Rankings are introduced, in order to be able to effectively represent
prime differential ideals by their characteristic sets, which allow
to test membership and solve other algorithmic problems. The latter
are defined as follows.

A differential polynomial $f$ is said to be reduced w.r.t. a differential 
polynomial $g$, if for any derivative operator $\theta\in\Theta$, the degree 
of $f$ w.r.t. to the the leader of $\theta g$ is less than that of $\theta g$.
A set of differential polynomials is called autoreduced, if each element
of this set is reduced w.r.t. the rest. Every autoreduced set is finite.
Among all autoreduced subsets of a given set of differential polynomials
$X$, choose one with the least set of ranks (which exists, because
$\le$ is a well-ordering on the set of finite sets of ranks).
This set is called a characteristic set of $X$. Note that the characteristic
set may not be unique, but its set of ranks is. 

Consider a differential polynomial $f$ as a polynomial in its leader $u$. Then
the leading coefficient of $f$ is called its initial, and the initial 
of any proper derivative $\theta f$, where $\theta\in\Theta\setminus\{1\}$,
is called the separant of $f$. For a finite set $\C$ of differential polynomials,
denote by $H_\C$ the product of the initials and separants of its elements.
It is a well-known fact (see e.g. \cite{Kol}) that every prime differential
ideal $P$ can be represented by its characteristic set $\C$ as
$$P=[\C]:H_\C^\infty,$$
that is, as the differential ideal generated by $\C$ and saturated by
the initials and separants of $\C$. We will extensively exploit this
representation in the paper.

\section{Main result}

\subsection{Different ways of posing the Ritt problem}
We give a list of equivalent formulations of the Ritt problem below and will use
one of them to attack the problem. We note that the equivalence takes place only 
for the fields $\K$ with a {\it splitting algorithm}, that is, an algorithm 
which, given a univariate polynomial $f\in \K[x]$, determines whether $f$ is
irreducible. This requirement is usually imposed in the context of polynomial
factorisation problems over $\K$ (see, e.g., \cite{Winkler, MillerGalois}). 
Examples of computable fields with splitting algorithm include the fields of 
rational numbers and rational functions, and their algebraic closures. An
example of a computable field without splitting algorithm is given in 
\cite{MillerGalois}: take a computably enumerable but non-computable
subset $S$ of the natural numbers, and consider the field $\Q[\sqrt{p_n}: n\in S]$, 
where $p_n$ is the $n$-th prime. By considering the field of fractions over
this field, we obtain a non-trivial differential field without splitting
algorithm. As in \cite[Section 5]{Miller08}, we suggest to think of the Ritt 
problem, especially of its third and fourth formulations in the theorem below, 
as a generalisation of the polynomial factorisation problem. 

\begin{theorem}
  The following problems are equivalent over a computable differential field
  of characteristic zero with a splitting algorithm (an algorithmic 
  solution for any one of them will also provide algorithms for the rest):
\begin{enumerate}
  \item Given a characteristic set of a prime differential 
    ideal, find a set of its generators.
  \item Given the characteristic sets of two prime differential 
    ideals $I_1$ and $I_2$, determine whether $I_1\subset I_2$.
  \item Compute a non-redundant prime decomposition of
    a radical differential ideal.
  \item Given a radical differential ideal specified by 
    a set of generators, determine whether it is prime. 
  \item Given a radical differential ideal specified by
    a set of generators, compute its prime decomposition together
    with the generators of the prime components as radical 
    differential ideals. 
  \item Given a radical differential ideal $I$ specified by
    a set of generators and a differential polynomial $f$, 
    determine whether $f$ is a zero-divisor modulo $I$. 
\end{enumerate}
\end{theorem}
\begin{proof}
  \begin{enumerate}
    \item[$1\Rightarrow 2$] Assume that we have an algorithm for
      finding generators of a prime differential ideal specified
      by its characteristic set. Applying this algorithm to the  
      characteristic set of $I_1$, compute its generator
      system $F_1$, i.e., $I_1=\{F_1\}$. Then $I_1\subset I_2$ 
      if and only if $F_1\subset I_2$. 
    \item[$2\Rightarrow 3$] Given a radical differential ideal
      $\{F\}$, one can apply the Ritt-Kolchin algorithm~\cite[Section IV.9]{Kol},
      \cite[Algorithm 5.5.15]{Pan}, \cite[Section 10]{Sit} to compute
      its prime decomposition $$\{F\}=P_1\cap\ldots\cap P_k,$$
      where each prime component is represented by its 
      characteristic set. 

      Note that the Ritt-Kolchin algorithm requires to test whether
      algebraic ideals of the form $(\C):H_\C^\infty$,
      where $\C$ is a coherent~\cite[Section III.8]{Kol}
      autoreduced set, are prime and, if not, to find two polynomials 
      not from the ideal, whose product
      belongs to the ideal. Apart from this test, the Ritt-Kolchin
      algorithm can be executed over any computable differential
      field of characteristic zero.

      The special case of $\C$ consisting of a single univariate
      polynomial shows that the existence of a splitting algorithm 
      over $\K$ is necessary for the above primality test. 
      We will show that it is also sufficient. For the proof
      of sufficiency of slightly stronger requirements see
      \cite[Section 5.5]{Pan}. 

      Indeed, if $\C$ is coherent autoreduced, then the ideal $(\C):H_\C^\infty$ 
      is radical \cite{Bou1,Fac}.
      A basis of this ideal can be computed over any computable field
      using Gr\"obner bases \cite[Section 4.4]{Cox}. Then, the generators of the 
      associated primes of the ideal can be computed over any computable 
      field with the splitting algorithm \cite{EHV}. If there is only one
      associated prime, the ideal is prime. Otherwise, by picking in each
      associated prime a polynomial that does not belong to the other associated
      primes (which can be done with Gr\"obner bases over any computable field),
      we can find the required product. 

      Assume that we have an algorithm for checking inclusion
      of prime differential ideals. Then in the above decomposition
      we can remove all redundant components, i.e., those $P_i$ 
      that contain another $P_j$. The resulting components will
      be essential, since any prime decomposition of a radical 
      ideal contains all essential components. Hence, we obtain 
      the non-redundant (essential) decomposition of $\{F\}$.
    \item[$3\Rightarrow 4$] Assume that we have an algorithm 
      for computing a non-redundant (=essential) prime decomposition
      of a radical differential ideal $\{F\}$. Then $\{F\}$ is
      prime if and only if this decomposition consists of one
      component.
    \item[$4\Rightarrow 1$] Given a characteristic set $\C$ of a prime 
      differential ideal $P$. Consider the algebraic ideals 
      $J_i=\left(\C^{(i)}\right):H_\C^\infty$, $i=0,1,2,\ldots$, where
      $$\C^{(i)}=\{\theta f\;|\;\ord\theta\le i,\;f\in \C\}.$$
      Let $F_i$ be a system of generators of $J_i$ (e.g., its
      Gr\"obner basis). By the basis theorem, there exists an index
      $i$ such that $\{F_i\}=P$.

      Assume that we have an algorithm for determining whether
      a radical differential ideal $\{F\}$ is prime. Applying
      this algorithm to $\{F_i\}$, $i=0,1,2,\ldots$, find 
      the least $i$ such that $\{F_i\}$ is prime. Then, since
      $\C$ is a characteristic set of $\{F_i\}$, we will have
      $\{F_i\}=P$, hence $F_i$ generates $P$ as a radical differential
      ideal.
    \item[$1\Rightarrow 5$] Given a radical differential ideal, one
      can apply the Ritt-Kolchin algorithm to compute its prime decomposition,
      where each prime component is represented by its characteristic
      set. Assuming that we have an algorithm for computing generators
      of a prime differential ideal, given its characteristic set,
      we obtain the generators of the components.
    \item[$5\Rightarrow 3$] Assume that we have an algorithm for
      computing a prime decomposition of a radical differential 
      ideal together with the generators of the prime components. 
      Then we can check inclusion of the components as in the proof
      $1\Rightarrow 2$ and thus exclude redundant components 
      as in $2\Rightarrow 3$.
    \item[$3\Rightarrow 6$] Assume that we have an algorithm for
      computing an essential prime decomposition $I=P_1\cap\ldots\cap P_k$
      of a radical differential ideal $I=\{F\}$. 

      \begin{lemma}
         A differential polynomial $f$ is a zero-divisor modulo $I$
         if and only if it belongs to an essential prime component of $I$.
      \end{lemma}
      \begin{proof}
      Let $f$ be a zero-divisor modulo $I$. Then by definition
      there exists a differential polynomial $g\not\in I$ such that 
      $fg\in I$. Since $g\not\in I$, there exists an essential prime
      component $P_i$ of $I$ which does not contain $g$, yet $fg\in I$
      implies $fg\in P_i$. Since $P_i$ is prime, we obtain $f\in P_i$.

      Vice versa, let $f$ be an element of an essential prime
      component of $I$. Let $P_{i_1},\ldots,P_{i_l}$ be the essential prime
      components of $I$ not containing $f$. Then for
      the ideal $I:f=\{g\;|\;gf\in I\}$, we have
      $$I:f=(P_1\cap\ldots\cap P_k):f=(P_1:f)\cap\ldots\cap(P_k:f)=P_{i_1}\cap\ldots\cap P_{i_l}.$$
      Since $P_1\cap\ldots\cap P_k$ is a non-redundant prime
      decomposition of $I$, we obtain $I:f\neq I$, hence $f$ is a
      zero-divisor modulo $I$. 
      \end{proof}
    \item[$6\Rightarrow 1$] The proof of this implication is similar
      to $[4\Rightarrow 1]$. 

      Given a characteristic set $\C$ of a prime 
      differential ideal $P$. Consider the algebraic ideals 
      $J_i=\left(\C^{(i)}\right):H_\C^\infty$, $i=0,1,2,\ldots$.
      Let $F_i$ be a system of generators of $J_i$.
      By the basis theorem, there exists an index
      $i$ such that $\{F_i\}=P$. Note that the product $H_\C$ 
      of initials and separants of $\C$ is not a zero-divisor 
      modulo $\{F_i\}$, since
      $$P=\{F_i\}\subseteq\{F_i\}:H_\C^\infty\subseteq [\C]:H_\C^\infty=P.$$

      Assume that we have an algorithm for determining whether
      a polynomial $f$ is a zero-divisor modulo a radical differential
      ideal $\{F\}$. Applying this algorithm to $H_\C$ and
      $\{F_i\}$, $i=0,1,2,\ldots$, find the smallest index $i$ such that the
      $H_\C$ is not a zero-divisor modulo $\{F_i\}$ (as we
      have shown above, such an index exists). 
      Then $\{F_i\}:H_\C^\infty=\{F_i\}$. On the other hand,
      $\{F_i\}:H_\C^\infty=P$, since $\C\subset \{F_i\}$. Thus,
      $\{F_i\}=P$. 
 \end{enumerate}
\end{proof}
\subsection{Canonical prime decomposition of a radical differential ideal}
In this section we will show how one can produce a prime decomposition of
a radical differential ideal $I$ given by a finite set of generators $F$ that does not
depend on this particular choice of generators. A good candidate for this decomposition
would be the essential prime decomposition but, again, it is still unknown whether 
the latter is computable. Our procedure is based on two simple observations
contained in Lemma~\ref{l1} and Proposition~\ref{p1}.

Assume that a ranking is fixed. Then a prime differential ideal has a canonical
characteristic set uniquely determined by the ideal \cite{Bou3,canonical,Pol,Dif}.
And given $F \subset \K\{Y\}$ we can compute some prime decomposition
\begin{equation}\label{decomposition}
I = \{F\} = \bigcap_{i=1}^k[\C_i]:H_{\C_i}^\infty,
\end{equation}
where $\C_i$ are the canonical characteristic sets of the corresponding prime
differential ideals.

\begin{lemma}\label{l1} Let $\C$ be a characteristic set of the highest rank among $\C_1,\ldots,\C_k$. Then the ideal $P := [\C]:H_\C^\infty$ is an essential prime component
of $I$ of the highest rank, among all essential prime components of $I$. 
\end{lemma}
\begin{proof}
Let $Q = [\B]:H_\B^\infty$, where $\B$ is the canonical characteristic set of $Q$, be an essential 
prime component of $I$ such that $P\supseteq Q$. Because of this inclusion and the definition of 
the characteristic set as an autoreduced subset of the least rank, we have $\rk\C \le \rk\B$.
Moreover, since $Q$ is an essential prime component of $I$, it must be among the ideals 
$[\C_i]:H_{\C_i}^\infty$, $i=1,\ldots,k$, which means that $\B=\C_l$ for some $l$, $1\Le l\Le k$.
Thus, due to the choice of $\C$, we have $\rk\C\ge\rk\B$. We conclude therefore that $\rk\C = \rk\B$.
Hence, according to \cite[Lemma 13]{canonical}, prime differential ideals $P$ and $Q$ are equal.
\end{proof}

\begin{proposition}\label{p1} Let $I=\{F\}$ be a radical differential ideal, and let 
$\C = C_1,\ldots,C_l$ be a finite set of differential polynomials satisfying
$\C\subset I\subset [\C]:H_\C^\infty$.
 Then{ }\footnote{Admitting
a slight abuse of notation in this formula, we denote by $F\cup H_\C$ 
the union of the set $F$ and the singleton set containing $H_\C$.}
$$
I = [\C]:H_\C^\infty\cap I:C_1\cap\ldots\cap I:C_l\cap\{F\cup H_\C\}.
$$
\end{proposition}
\begin{proof} Let 
$$f \in [\C]:H_\C^\infty\cap I:C_1\cap\ldots\cap I:C_l\cap\{F\cup H_\C\}.$$ Then 
$$C_i\cdot f \in I,\ 1\Le i\Le l$$ and there exists $h \in H_\C^\infty$ such that
$$h\cdot f\in [\C].$$ 
We then have 
$$
h\cdot f^2 \in f\cdot[\C] \subset \{f\cdot C_1,\ldots,f\cdot C_l\} \subset I.
$$
and, therefore, $f \in \{F\}:H_\C^\infty$.
Since due to, for example, \cite[Proposition 2.1]{Fac}
$$\{F\} = \{F\}:H_\C^\infty\cap\{F\cup H_\C\},
$$ we conclude that
$f \in I$. The other inclusion follows from $I \subset [\C]:H_\C^\infty$.
\end{proof}

If in decomposition~\eqref{decomposition} there are several characteristic sets of the 
highest rank, say, $\C_1,\ldots,\C_q$, and $\C_i = C_{i,1},\ldots,C_{i,p_i}$, $1\Le i \Le q$, then
\begin{align}\label{decomposition2}
I =[\C_1]:H_{\C_1}^\infty&\cap I:C_{1,1}\cap\ldots\cap I:C_{1,p_1}\cap\left\{F\cup H_{\C_1}\right\}\cap\ldots\\
&\cap \left[\C_q\right]:H_{\C_q}^\infty\cap I:C_{q,1}\cap\ldots\cap I:C_{q,p_q}\cap\left\{F\cup H_{\C_q}\right\}.\notag
\end{align}
Note that $\C_1,\ldots,\C_q$ and, therefore, all ideals in~\eqref{decomposition2} are uniquely determined by $I$, that is, they do not depend on the choice of generators
of $I$. Moreover, the ideals 
\begin{align}\label{id1}
I:C_{i,j}
\end{align} 
and
\begin{align}\label{id2}
\left\{F\cup H_{\C_i}\right\}
\end{align}
strictly
contain $I$, $1\Le i \Le q$, $1\Le j \Le p_i$. Indeed, all $C_{i,j}$ are elements of essential
prime components of $I$, therefore $I:C_{i,j}\supsetneq I$. And, since $H_{\C_i}$ does not belong to
the corresponding essential prime component $[\C_i]:H_{\C_i}^\infty$, it does not belong
to $I$, and we have $\{F\cup H_{\C_i}\}\supsetneq I$.

By the Ritt-Raudenbush theorem,
a strictly increasing chain of radical differential ideals terminates. Therefore, by
computing recursively the canonical, generator-independent prime decomposition
of~\eqref{id1} and~\eqref{id2}, we obtain a generator-independent decomposition
of the original ideal $I$. We note that the prime decomposition required 
in step~\eqref{decomposition} can be computed for the radical ideals given in 
a saturated form, as in~\eqref{id1}.

\subsection{Ranking-independent canonical decomposition and generators}

Even though the canonical prime decomposition computed by the method described in the previous 
section is independent of the given generators of the radical ideal $I$, it does 
depend on the choice of ranking. It is possible to obtain a ranking-independent canonical
prime decomposition via the following modification. As in~\eqref{decomposition}, compute
any prime decomposition of $I$. Then, instead of extracting the characteristic sets of the 
highest rank w.r.t. a fixed chosen ranking (this step is ranking-dependent), take
all prime components whose canonical characteristic sets have the highest rank w.r.t. {\it some} 
ranking, and consider all {\it these} characteristic sets instead of the above $\C_1,\ldots,\C_q$.
This step can be accomplished by computing the universal characteristic set \cite{ucs} for each 
prime component. The remaining steps are the same.

Finally, from the canonical prime decomposition 
 $$I = \bigcap_{i=1}^k[\C_i]:H_{\C_i}^\infty$$
one can obtain a canonical set of generators of the radical differential ideal, which
depends only on this ideal. For $j=0,1,2,\ldots$ consider the algebraic ideal
  $$I_j = \bigcap_{i=1}^k\left(\C_i^{(j)}\right):H_{\C_i}^\infty$$
and compute its Gr\"obner basis $\B_j$ w.r.t. the lexicographic term order induced 
on power products of derivatives by the ranking (if a ranking-independent
canonical set of generators is sought, compute the universal Gr\"obner basis).
Stop at the least value of $j$ such that $\{\B_j\}=I$, and output $\B_j$. Such $j$ 
exists, because $I=\bigcup_{j=0}^\infty I_j$. Note that the
equality of two radical differential ideals given by generators can
be checked by testing membership of generators of one ideal to the other ideal.

\section{Acknowledgemenets}
We thank Michael Singer and William Sit for their helpful suggestions and support, 
and Russell Miller for the introduction to the theory of computable fields and 
subsequent discussions of this topic in connection with the Ritt problem.

\bibliographystyle{amsplain}
\bibliography{rittprob}
  
\end{document}